\documentclass[10pt]{amsart}
\usepackage{amsmath,amsthm,amssymb,bbm,eufrak,here}

\theoremstyle{plain}
\newtheorem{thm}{Theorem}[section]

\newtheorem{cor}[thm]{Corollary}
\newtheorem{prop}[thm]{Proposition}
\newtheorem{lemma}[thm]{Lemma}

\renewcommand\mathbb[1]{\mathbbm{#1}}
\newcommand\V{\ensuremath{\mathbb{V}}}
\renewcommand\P{\ensuremath{\mathbb{P}}}
\newcommand\Z{\ensuremath{\mathbb{Z}}}
\newcommand\I{\ensuremath{\mathcal{I}}}

\newcommand\K{\ensuremath{K_W}}

\newcommand\Span[1]{span\langle #1 \rangle}
\newcommand\corresponds{\rightsquigarrow}

\newcommand\comment[1]{}

\begin{document}

\title{Hessian quartic surfaces that are Kummer surfaces}
\author{Joel E. Rosenberg}

\begin{abstract}
In 1899, Hutchinson \cite{hutch} presented a way to obtain a
three-parameter family of Hessians of cubic surfaces as blowups of
Kummer surfaces. We show that this family consists of those Hessians
containing an extra class of conic curves. Based on this, we find the
invariant of a cubic surface $C$ in pentahedral form that vanishes if
its Hessian is in Hutchinson's family, and we give an explicit map
between cubic surfaces in pentahedral form and blowups of Kummer
surfaces.
\end{abstract}

\maketitle

\section{Introduction}

Let $C$ be a cubic surface in $\P^3$. Among the many interesting
geometrical objects associated to $C$ is its Hessian, a quartic
surface $H$ in $\P^3$. It was found in the nineteenth century
\cite{segre} that $H$ will have ten double points, and will contain
ten lines through those points. Conversely, it was shown that any
irreducible quartic surface containing an appropriate configuration of
lines and double points would be the Hessian of a unique cubic
surface.

Another class of objects of interest in classical algebraic geometry
is the class of Kummer surfaces. Given an abelian surface $A$, we have
an action of the group $\{\pm1\}$ by multiplication, and we can take
the quotient $K=A/\{\pm1\}$, a surface with 16 double points. On a
Kummer surface, there are 16 special curves, called tropes, each of
which passing through six of the double points. It was found
\cite{hutch} that one can choose a certain other subset $W$ of 6 of
the double points, called a Weber hexad, and blow these points up, to
obtain a surface \K\ that embeds in $\P^3$ as the Hessian of a cubic
surface. That is, there remain 10 double points on \K, and one can
embed the surface such that of the sixteen tropes, ten are taken to
straight lines in $\P^3$, and in the same configuration as referred to
above. So we conclude that \K\ is the Hessian of some cubic surface.

Now, it is known that the moduli space of cubic surfaces is a
four-dimensional normal variety, while there is only a three-parameter
family of Kummer surfaces, each of which having only finitely many
Weber hexads. So we can hope that the locus of cubic surfaces whose
Hessians are Kummer in the above sense will be a divisor in the space
of all cubic surfaces. In this paper we prove

\begin{thm}\label{invt}
Let $\P^3$ be taken to be the hyperplane in $\P^4$ with $\sum_{i=0}^4
X_i=0$. Let $C\subset\P^3$ be the cubic surface
$\V(\sum_{i=0}^4\frac1{\mu_i}X_i^3)$, where
$\mu_0\mu_1\mu_2\mu_3\mu_4\not=0$, and assume $C$ is smooth. Then the
Hessian of $C$,
$$\V(\mu_0X_1X_2X_3X_4+\mu_1X_0X_2X_3X_4+\mu_2X_0X_1X_3X_4+
\mu_3X_0X_1X_2X_4+\mu_4X_0X_1X_2X_3),$$
is the blowup of a Weber hexad on a Kummer surface if and only if the
coefficients $\mu_i$ satisfy the following irreducible cubic
condition:
$$\sum_{i=0}^4\mu_i^3-\sum_{i\not=j}\mu_i^2\mu_j
+2\sum_{i\not=j\not=k}\mu_i\mu_j\mu_k=0.$$
\end{thm}

The reader familiar with the invariant theory of cubic surfaces will
be curious how this locus fits in with the classical invariants. We
can interpret the previous result in that language, obtaining the
following.
\begin{cor}\label{thirtytwo}
Let $\P^{19}$ be the parameter space of cubic forms on $\P^3$. Let $X$
be the locus in $\P^{19}$ of cubic surfaces whose Hessians are
isomorphic to blowups of Weber hexads on Kummer surfaces, embedded as
in Theorem~\ref{invt}. Then $X$ is $SL(4)$-invariant, and the closure
of $X$ is a divisor in $\P^{19}$. If we label the classical invariants
as $I_8,I_{16},I_{24},I_{32},I_{40}$, following \cite{hunt}, then the
polynomial on $\P^{19}$ given by
$$I_8I_{24}+8I_{32}$$
is irreducible, is degree 32, and vanishes on $X$.
\end{cor}

This result, while somewhat satisfying, is only the beginning of the
story. The next natural question to ask is, ``If this divisor in the
parameter space of $\mu_i$s is associated to the moduli space of
Kummer surfaces with Weber hexads, then what is this correspondence?''
We also answer this question. Recall that the generic abelian surface
is the Jacobian of a unique genus 2 curve, say $A=J(B)$, and that $B$
can be specified as the double cover of $\P^1$ branched at six points
$a,b,c,d,e,f$. It turns out that upon placing an ordering on these six
points, we can specify a unique Weber hexad in $K$. We then have the
following theorem.

\begin{thm}
Let $a,b,c=0,d=1,e,f=\infty$ be six distinct points in $\P^1$. Let $B$
be the double cover of $\P^1$ branched at these six points, and let
$A$ be the abelian surface $J(B)$. Let $W\subset A$ be the Weber hexad
$$\{0,b+c-2a,c+d-2a,d+e-2a,e+f-2a,f+b-2a\}.$$
Then the surface \K\ obtained by blowing up $K=A/\{\pm1\}$ at $W$ can
be embedded in $\P^3$ as the Hessian of the surface
$$\V(\sum_{i=0}^4\frac1{\mu_i}X_i^3),$$
where the coefficients $\mu_i$ are given by
\begin{align*}
\mu_0&=a(1-b),&
\mu_1&=e(1-a),&
\mu_2&=b(e-a),&
\mu_3&=(e-b),&
\mu_4&=(a-b)(1-e).
\end{align*}

Conversely, if $\mu_0\mu_1\mu_2\mu_3\mu_4\not=0$ and
$$\sum_{i=0}^4\mu_i^3-\sum_{i\not=j}\mu_i^2\mu_j
+2\sum_{i\not=j\not=k}\mu_i\mu_j\mu_k=0,$$
then let
\begin{align*}
a&=\frac{\mu_0+\mu_3+\mu_4-\mu_1-\mu_2}{2\mu_3},\\
b&=\frac{2\mu_2}{\mu_1+\mu_2+\mu_3-\mu_0-\mu_4},\\
e&=\frac{\mu_0+\mu_3-\mu_4-\mu_1-\mu_2}{\mu_0+\mu_3+\mu_4-\mu_1-\mu_2}
\end{align*}
be points in $\P^1$. If these points are all distinct, and none of
these points are $0$, $1$, or $\infty$, then let $B$ be the double
cover of $\P^1$ branched at $a,b,c=0,d=1,e,f=\infty$, and let \K\ be
the blown up Kummer specified above. Then the Hessian surface
$$H=\V(\mu_0X_1X_2X_3X_4+\mu_1X_0X_2X_3X_4+\mu_2X_0X_1X_3X_4+
\mu_3X_0X_1X_2X_4+\mu_4X_0X_1X_2X_3)$$
is isomorphic to \K.
\end{thm}

This still is not the whole story as it should be. If one recalls that
any smooth cubic surface can be obtained as the blow up of $\P^2$ at
six points, the question arises, given a genus 2 curve $B$ and a Weber
hexad $W\subset J(B)$, what six points we should blow up in $\P^2$ to
obtain a cubic surface $C$ whose Hessian $H$ is isomorphic to \K. As
far as this author is aware, this question has not been answered.
Failing that, we present proofs of the above two theorems, in the hope
that our techniques can be extended to answer the remaining questions
about Kummer Hessian surfaces.

\section{The geometry of Hessian quartic surfaces}\label{hessian}

We will begin our exploration by collecting some results about the
geometries of Hessian surfaces and Kummer surfaces, and then apply
them to the theorems at hand. First, let us describe some of the
geometry of the Hessian of a generic cubic surface. Let $\P^3$ be the
hyperplane
$$\V(\sum_{i=0}^4 X_i)\subset\P^4,$$
as above, and for $0\leq i<j\leq4$, let
$$\ell_{ij}=\V(X_i,X_j)\subset\P^3.$$
Let $L=\bigcup_{ij}\ell_{ij}$. Then we have the following lemma.
\begin{lemma}
Let $H$ be a quartic form on $\P^3$ that vanishes on $L$. Then $H$ is
in the linear span of
$$\langle X_0X_1X_2X_3,X_0X_1X_2X_4,X_0X_1X_3X_4,X_0X_2X_3X_4,X_1X_2X_3X_4\rangle.$$
As a result, $H$ is double at the 10 points $p_{ijk}=\V(X_i,X_j,X_k)$.
\end{lemma}

\begin{proof}
For the first statement, we observe that $\V(H,X_0)$ by assumption
contains the four lines $\ell_{01},\ell_{02},\ell_{03},\ell_{04}$. So
$H$ lies in the ideal $(X_0,X_1X_2X_3X_4)$. Applying this
symmetrically yields the result. The second statement is immediate.
\end{proof}

Observe that if $\mu_0\mu_1\mu_2\mu_3\mu_4\not=0$, then the cubic
surface $\V(\sum_{i=0}^4\frac1{\mu_i}X_i^3)$ has Hessian equal to 
$$H=\mu_0X_1X_2X_3X_4+\mu_1X_0X_2X_3X_4+\mu_2X_0X_1X_3X_4+
\mu_3X_0X_1X_2X_4+\mu_4X_0X_1X_2X_3.$$
Clebsch showed \cite{salmon} that the generic cubic surface is
isomorphic to an essentially unique cubic surface of this
(pentahedral) form, but this result is largely irrelevant to our work
here, so simply recall the result, and do not pursue it further.

So a Hessian in this family contains ten lines and ten double points,
such that each line passes through three of the double points and each
double point lies on three of the lines. We now give one result each
about the double points and the lines of a Hessian quartic $H$.

\begin{prop}
Projection away from the node $p_{012}=\V(X_0,X_1,X_2)$ gives a
rational map from $\V(H)$ to $\P^2$, which is generically 2-to-1. The
branch locus of this map is the union of two cubic curves in $\P^2$,
with equations
\begin{multline*}
\left(sX_0X_1X_2+(\mu_0X_1X_2+\mu_1X_0X_2+\mu_2X_0X_1)(-X_0-X_1-X_2)\right)\\
\left(\bar{s}X_0X_1X_2+(\mu_0X_1X_2+\mu_1X_0X_2+\mu_2X_0X_1)(-X_0-X_1-X_2)\right)
\end{multline*}
where $s$ and $\bar{s}$ are the roots of
$$s^2-2(\mu_3+\mu_4)s+(\mu_3-\mu_4)^2=0.$$
\end{prop}
\begin{proof}
The first statement follows because a generic line through $p_{012}$
meets $\V(H)$ twice there, and at two other points. The second
statement follows by computing the discriminant of $H$, viewed as a
quadric in $X_3$.
\end{proof}

Observe that the two cubics given are tangent to the quadric
$\mu_0X_1X_2+\mu_1X_0X_2+\mu_2X_0X_1=0$ at the three points
$\V(X_i,X_j)_{0\leq i<j\leq2}$, and so meet each other to order two
there. These points are the images of the lines in $H$ through
$p_{012}$. The cubic curves have their remaining three intersections
transverse, all along the line $X_0+X_1+X_2=0$. This line is the image
of the line $\ell_{34}$. We next prove a result about this line.

\begin{prop}
The plane $\mu_4X_3+\mu_3X_4=0$ is tangent to $H$ at every point of
the line $\ell_{34}$.
\end{prop}
\begin{proof}
This is immediate from the equation of the surface.
\end{proof}

Observe that the intersection of this plane with the surface then
consists of the line $\ell_{34}$, counted twice, and a conic with
equation $\mu_0X_1X_2+\mu_1X_0X_2+\mu_2X_0X_1=0$, the same conic
referred to above.

Finally, observe that if a cubic surface $C$ has a node at a point
$p$, then its Hessian is also nodal at $p$, with the same tangent
cone. Connversely, if a Hessian in our four-parameter family acquires
a node other than the ten coordinate points, the corresponding cubic
surface also acquires a node. So, if we restrict our attention to
Hessians $H$ of smooth cubic surfaces, we may assume that $H$ contains
only ten nodes, and that the discriminant sextics described above are
smooth away from the six images of the nodes.

\section{The geometry of Kummer surfaces}

For this section, we will largely follow the development in \cite{gh},
with one exception. Since we have already made use of subscripted
numbers for our cubic surface in pentahedral form, we will begin with
six distinct points labelled $a,b,c,d,e,f\in\P^1$. So begin with these
points, and let $B\to\P^1$ be the genus 2 curve that is the double
cover branched over these six points. We will also label the
ramification points in $B$ by the letters $a\ldots f$. Then the
Jacobian of $B$ is an abelian surface $A$, with 16 two-torsion points,
and these correspond to the divisors
$$0,\qquad\{b-a,c-a,\dots,f-a\},\qquad\{b+c-2a,\dots,e+f-2a\}.$$

Recall that a theta divisor on $A$ is an image of the curve $B$ under
a map $p\mapsto p-D$, where $D$ is some divisor of degree 1 on $B$. If
$D$ is any ``two-torsion'' point, i.e., if $2D\sim2a$, then the
theta-divisor given will pass through 6 of the two-torsion points of
$A$. We will refer to these 16 divisors on $A$ as tropes, and will
label the trope corresponding to the divisor $D$ by the symbol
$\Theta_D$. Note
that they give 16 distinguished subsets of 6 two-torsion points. We
will now define a different sort of set of six two-torsion points,
called a Weber hexad. A Weber hexad is a set of 6 points of the
two-torsion of $A$ such that 10 of the tropes each contain 3 of the
points of the hexad, and the other 6 tropes each contain exactly one
of the points of the hexad. For example, the six points
$$0,b+c-2a,c+d-2a,d+e-2a,e+f-2a,f+b-2a$$
have this property: only 0 lies on $\Theta_a$, etc. On a given abelian
surface, there are exactly 192 Weber hexads, which can be obtained
from the one above (each one 60 times) by acting by translation by the
two-torsion of $A$, and by acting on $a\ldots f$ with the group $S_6$.

Now, as discussed in the introduction, if we identify points on $A$
with their negatives, we obtain a surface $K$ with sixteen double
points, the image of the two-torsion. We can desingularize these nodes
by blowing up, or equivalently by blowing up the two-torsion points on
$A$ before taking the quotient. Observe that if we blow up the six
points of a Weber hexad, we will be left with a surface with 10 nodes,
just like our Hessians above. In fact, one makes the following claim,
first noticed by Hutchinson \cite{hutch}.

\begin{prop}
Let $a\ldots f, B, A, K$ be as above. Let
$$W=\{0,b+c-2a,c+d-2a,d+e-2a,e+f-2a,f+b-2a\}\subset A.$$
If one maps $A$ to projective space using the linear series
$|4\Theta_a-2W|$, one gets a map to a quartic surface $\K\subset\P^3$,
with 10 nodes, with the following properties. There exist 5 planes in
$\P^3$ such that \K\ is nodal at the intersection of any three of
these planes, and contains the line that is the intersection of any
two of these planes. Further, the image of each of the point of $W$ is
a conic in \K.
\end{prop}

\begin{proof}
One may check using homological criteria that the linear series
$|4\Theta_a-2W|$ has rank 4, and that all of its sections are even
functions, so identify points with their negatives. So, the map is
$\rho:A\dashrightarrow\P^3$, and factors through $K$. The points of
$W$ are in the base locus, so get blown up. Now, the self-intersection
$(4\Theta_a-2W)^2=8$, which divided by two gives 4, so \K\ is a
quartic surface, as stated. Also, if $E_0$ is the exceptional divisor
over $0\in A$, then $E_0$ is part of the ramification locus of the
map. So if $C_0$ is the image of $E_0$ in $\P^3$, and $\omega$ the
hyperplane class in $\P^3$, then
$C_0\cdot\omega=E_0\cdot(4\Theta_a-2W)=2$, and $C_0$ is a conic.
Likewise the remaining points of $W$ also map to conics.

Now, for clarity we will introduce the notation $p_\lambda$ for the
two-torsion point we have been calling $\lambda\in A$. 
Observe that 
\begin{align*}
&(\Theta_b-p_0-p_{b+c-2a}-p_{e+b-2a})+\\
&(\Theta_d-p_0-p_{d+e-2a}-p_{c+d-2a})+\\
&(\Theta_{b+c-a}-p_{b+c-2a}-p_{d+e-2a}-p_{c+d-2a})+\\
&(\Theta_{c+d-a}-p_{c+d-2a}-p_{e+a-2a}-p_{d+e-2a})
=4\Theta_a-2W.
\end{align*}
So these four tropes have coplanar image, and are all lines. Letting
the group $\Z/5\Z$ act by $(bcdef)$, we get five such planes, and the
result.
\end{proof}

Hutchinson largely ignores the conics coming from $W$, because his
purpose is to study low-degree curves that arise on every Hessian
quartic, not solely on the Kummer surfaces. In the next section, we
will take the opposite approach, and study those Hessian surfaces that
do contain conics like these, and find that this extra class of curves
is enough to make a surface Kummer.

\subsection{Labelling}\label{label}

Before we proceed to our main results, we pause here for a discussion
of labelling. To conform with our names for the 5 planes in the
discussion of the Hessian, we will let
\begin{align*}
P_0=\V(X_0)&=\Span{\rho(p_b),\rho(p_c),\rho(p_d)},\\
P_1=\V(X_1)&=\Span{\rho(p_c),\rho(p_d),\rho(p_e)},\\
P_2=\V(X_2)&=\Span{\rho(p_d),\rho(p_e),\rho(p_f)},\\
P_3=\V(X_3)&=\Span{\rho(p_e),\rho(p_f),\rho(p_b)},\\
P_4=\V(X_4)&=\Span{\rho(p_f),\rho(p_b),\rho(p_c)}.
\end{align*}
Then we may label the lines on \K\ as $\ell_{ij}$ as before, and
obtain another set of names for the ten nodes. However, the real
interest lies in what we can say about the conics. Observe that the
conic $C_0$ meets the tropes
$$\Theta_d,\Theta_f,\Theta_c,\Theta_e,\Theta_b,$$ that is, the lines
$$\ell_{02},\ell_{24},\ell_{41},\ell_{13},\ell_{30}.$$
So, we would do well to associate to $C_0$ the \textit{cyclic
ordering} $(02413)$. Similarly, one finds that the other exceptional
divisors should be assigned the labels
\begin{align*}
C_{b+c-a}&\corresponds(03214),\\
C_{c+d-a}&\corresponds(01432),\\
C_{d+e-a}&\corresponds(04312),\\
C_{e+f-a}&\corresponds(01324),\\
C_{f+b-a}&\corresponds(03421).
\end{align*}
This accounts for 6 of the 12 cyclic orders on five letters. The
astute reader will ask to what do the other cyclic orders correspond.
The astute reader will also have noticed that a conic is a plane
curve, and since a plane meets a quartic surface in a degree 4 curve,
each conic on \K\ must be coplanar with a residual conic on \K. The
answer, of course, is that we should label these residual conics with
the complementary orderings.

Now take one of these conics, say the one labelled $(ijklm)$, and
consider how it meets the lines on \K\ through a node $p_{rst}$. If
$rst$ are consecutive letters in the cyclic order $(ijklm)$, then the
conic will meet the lines $\ell_{rs}$ and $\ell_{st}$, but not
$\ell_{rt}$. If the letters $rst$ are not consecutive in $(ijklm)$,
for example, $ijl$, then the conic will only meet the one line
$\ell_{ij}$, while its residual will meet the other two lines.

In summary, the plane containing $C_{(ijklm)}$ and $C_{(ikmjl)}$
corresponds to one of the 6 subgroups of order 5 in $S_5$. The
elements of order 5 are all conjugate under the action of $S_5$, but
break into two orbits under the action of $A_5$, with each element
conjugate to its inverse. In one of these $A_5$-orbits we get the 6
cyclic orders we assigned to the exceptional divisors in \K, while in
the other $A_5$-orbit we get the labels of the residual conics.

\section{Conics on Hessian quartic surfaces}

Now, we return to the question of when a Hessian quartic surface
contains ``extra'' conics. We observed in section \ref{hessian} that
every Hessian quartic surface contains 10 conics, each lying in the
tangent plane to one of the ten lines. We saw that the conic in the
plane of $\ell_{34}$ met all three of the lines through the node
$p_{012}$. In the last section, we saw that in a Kummer Hessian
surface, there existed twelve other conics, not coplanar with any
node, and that for each node, each of these conics met exactly one or
two of the lines through that node. So, we assume that we have a
surface $H$ in the four-parameter family we have been studying, that
$H$ contains a conic $C$ in its smooth locus such that $C$ meets the
lines $\ell_{01}$ and $\ell_{02}$, but not $\ell_{12}$. We will find
that this condition is necessary and sufficient for the Hessian to be
Kummer, and so prove Theorem~\ref{invt}. We begin with the following
lemma.

\begin{lemma}
Let $H$ be the quartic surface
$$\V(\sum_{i=0}^4\mu_i\prod_{j\not=i}X_j),$$
and assume $\prod_{i=0}^4\mu_i\not=0$, so $H$ is a Hessian, and that
$H$ has only the ten nodes it should, i.e., the cubic surface is
smooth. Assume there exists a conic $C$ contained in the smooth locus
of $H$ such that $C$ meets the lines $\ell_{01}$ and $\ell_{02}$, but not
$\ell_{12}$. Then the $\mu_i$ satisfy the following irreducible cubic
form:
$$\sum_{i=0}^4\mu_i^3-\sum_{i\not=j}\mu_i^2\mu_j
+2\sum_{i\not=j\not=k}\mu_i\mu_j\mu_k=0.$$
\end{lemma}

\begin{proof}
To begin with, since $C$ meets $\ell_{01}$ at a smooth point, it must
be tangent to the plane $\mu_1X_0+\mu_0X_1=0$, and likewise tangent to
the plane $\mu_2X_0+\mu_0X_2=0$. Now, let $\pi:\K\to\P^2$ be the
projection away from the point $p_{012}$, and let $Q=\pi(C)$, a conic.
Our observations have placed four linear conditions on $Q$, so $Q$
must lie in the pencil
$$(\mu_0X_1X_2+\mu_1X_0X_2+\mu_2X_0X_1)+\alpha X_0^2.$$
Now our goal is to find which $Q$ in this pencil can have a conic in
its preimage. Since $Q$ passes through the images of the lines
$\ell_{01}$ and $\ell_{02}$, and with the right tangent direction, the
pullback $\pi^{-1}(Q)$ will always contain these lines, each counted
twice. The rest of the preimage will then be a quartic curve,
dominating $Q$ and mapping 2-to-1 to it. This quartic will be branched
over $Q$ at four points, the remaining intersections of $Q$ with the
discriminant locus away from the lines through $p_{012}$. To have the
preimage decompose, these four intersections must coincide in pairs.

So for each of the cubic curves making up the branch locus, we ask
which elements of the pencil meet it a third time non-reducedly. So
let 
$$E_s=\V\left(sX_0X_1X_2+(\mu_0X_1X_2+\mu_1X_0X_2+\mu_2X_0X_1)(-X_0-X_1-X_2)\right).$$
Since $H$ is assumed to have only the ten nodes, $E_s$ will be smooth,
so an elliptic 
curve, and there will be exactly four elements of the pencil where the
$g^1_2$ given by residual intersection with the conic branches. Two of
these are uninteresting: we already know the preimage of the conic
$$\mu_0X_1X_2+\mu_1X_0X_2+\mu_2X_0X_1,$$
and the double line $X_0^2$ pulls back to the plane $P_0$. So we are
left with a quadratic equation in $\alpha$ indicating which conics
meet $E_s$ interestingly. This quadratic is
$$T(s,\alpha)=4\mu_0\alpha^2
+[(s-\mu_0)^2+(\mu_2-\mu_1)^2-2(s+\mu_0)(\mu_2+\mu_1)]\alpha
+4s\mu_1\mu_2.$$

We want for there to be a conic that meets both $E_s$ and
$E_{\bar{s}}$ interestingly, so we take the resultant of $T(s,\alpha)$
and $T(\bar{s},\alpha)$, to find for what values of
$\mu_0,\mu_1,\mu_2,\mu_3,\mu_4$ these quadratics have a common
solution. We find that the resultant is
$$512\mu_0\mu_1\mu_2\mu_3\mu_4\left[
\sum_{i=0}^4\mu_i^3-\sum_{i\not=j}\mu_i^2\mu_j
+2\sum_{i\not=j\not=k}\mu_i\mu_j\mu_k
\right],$$
and since we are not interested in the cases where some $\mu_i=0$, we
keep only the last factor.  
\end{proof}

Observe that this form is symmetric in the five variables. Now, since
this cuts out an irreducible threefold in the space of $\mu$s, we have
almost proven Theorem~\ref{invt}. We deal with the remaining issues
below. But while we have the quadratics $T(s,\alpha)$ and
$T(\bar{s},\alpha)$ in hand, we observe that their difference is
linear in $\alpha$, and solving, we may write
$$\alpha=\frac{2\mu_1\mu_2}{\mu_0+\mu_1+\mu_2-\mu_3-\mu_4}.$$

We now return to the theorem, which we restate as follows:
\begin{thm}
Let $H$ be the quartic surface
$$\V(\sum_{i=0}^4\mu_i\prod_{j\not=i}X_j),$$
and assume $\prod_{i=0}^4\mu_i\not=0$, so $H$ is a Hessian, and that
$H$ has only the ten nodes it should, i.e., the cubic surface is
smooth. Assume the $\mu_i$ satisfy the cubic form
$$\sum_{i=0}^4\mu_i^3-\sum_{i\not=j}\mu_i^2\mu_j
+2\sum_{i\not=j\not=k}\mu_i\mu_j\mu_k=0.$$
Then there exists a Kummer surface $K$ and a Weber hexad $W$ such that
$H\cong\K$.
\end{thm}

\begin{proof}
We begin by showing that for any $H$ satisfying this hypothesis, there
are only finitely many conics satisfying the hypothesis of the lemma.
But this is easy, since at most two choices of $\alpha$ can give
acceptable image conics $Q$, and each of these will have at most two
preimages in $H$. In fact, we observe that since $s\not=\bar{s}$, only
one choice of $\alpha$ will work.

Next, we show that from the existence of such a conic $C$, we can
deduce the existence of eleven others meeting different subsets of the
lines, as in section~\ref{label}. We know that $C$ meets the lines
$\ell_{01}$ and $\ell_{02}$. If we look at its intersection with the
plane $P_2$, we see that it must hit another line, which we will
assume is $\ell_{23}$. We then intersect it with the plane $P_4$, and
conclude that it must hit $\ell_{14}$ and $\ell_{34}$. So we label our
conic $C_{(01432)}$. As in the proof of the lemma, we may project away
from $p_{012}$ and pull back, residuating $C_{(01432)}$ to a conic,
which we will call $C_{(01342)}$, and we can verify that it meets the
appropriate lines and so deserves that name. Similarly, we can project
the Hessian away from $p_{014}$, $p_{134}$, $p_{234}$, and $p_{023}$,
and obtain four other conics. Residuating each of these six conics in
the plane containing them gives us our total of twelve. Using the
observation that only two conics will meet, $\ell_{01}$ and
$\ell_{02}$, but not $\ell_{12}$, and the image under $S_5$ of this
fact, we can conclude that these are the only twelve interesting
conics on $H$.

Now to prove the theorem, we observe that the conics in one
$A_5$-orbit, say
$$C_{(02413)},C_{(03214)},C_{(01432)},C_{(04312)},C_{(01324)},C_{(03421)},$$
are all disjoint, and each had self-intersection $-2$. So we may blow
these down to obtain a 16-nodal K3 surface, which is known to be Kummer.
\end{proof}

Given this theorem, we may now restate, and prove,
Corollary~\ref{thirtytwo}.
\begin{cor}
Let $\P^{19}$ be the parameter space of cubic forms on $\P^3$. Let $X$
be the locus in $\P^{19}$ of cubic surfaces whose Hessians are
isomorphic to blowups of Weber hexads on Kummer surfaces, embedded as
in Theorem~\ref{invt}. Then $X$ is $SL(4)$-invariant, and the closure
of $X$ is a divisor in $\P^{19}$. If we label the classical invariants
as $I_8,I_{16},I_{24},I_{32},I_{40}$, following \cite{hunt}, then the
polynomial on $\P^{19}$ given by
$$I_8I_{24}+8I_{32}$$
is irreducible, is degree 32, and vanishes on $X$.
\end{cor}

\begin{proof}
Inside $\P^{19}$, consider the 4-plane $P$ of cubic forms
$$\sum_{i=0}^4\lambda_iX_i^3,$$
where as usual $X_0+X_1+X_2+X_3+X_4=0$. We have been concentrating
our attention on the locus in this family where no $\lambda_i$ equals
zero, and have been writing $\mu_i=\frac1{\lambda_i}$.
If we pull back our cubic condition
$$\sum_{i=0}^4\mu_i^3-\sum_{i\not=j}\mu_i^2\mu_j
+2\sum_{i\not=j\not=k}\mu_i\mu_j\mu_k=0$$
in the $\mu_i$ variables to the 4-plane $P$, we see that we have
described an open subset of a degree 12 threefold
$T\subset P\subset\P^{19}$. The locus $X$ is the $SL(4)$-orbit of $T$,
and since $T$ has finite stabilizer in $SL(4)$, the closure of $X$ is
a divisor in $\P^{19}$.

We next look for those $SL(4)$-invariants on $\P^{19}$ which vanish on
$X$. These must vanish on $T$, a degree 12 threefold inside a $\P^4$
in $\P^{19}$. The ring of invariants of cubic forms has no element of
degree 12, but if we refer to \cite{hunt} or \cite{salmon} to recall
how the invariants restrict to $P$, we will
find that the invariant
$$\I=I_8I_{24}+8I_{32}$$
vanishes on $T$. Indeed, this irreducible degree 32 invariant cuts out
the closure of $X$ in $\P^{19}$, and \I\ restricts on $P$
to an irreducible degree 12 polynomial, multiplied by
$(\lambda_0\lambda_1\lambda_2\lambda_3\lambda_4)^4$.
\end{proof}

\section{Finding the correspondence}

So given a Hessian quartic surface satisfying the hypothesis of
Theorem~\ref{invt}, we know it to be \K\ for some
choice of Kummer surface and Weber hexad, and the next question is
to which genus 2 curve it corresponds and to which Weber hexad. We
answer this by observing that the lines $\ell_{ij}$ correspond to
tropes on the Kummer, and that the three nodes on a line and the right
choice of three places where conics meet the line give six points on
$\P^1$ which specify the genus 2 curve $B$. So our task is to
explicitly find the six planes containing the conics.

To do this, we return to the pullback of the conic $Q$ of the last
section. So for any Hessian $H$, not necessarily Kummer, let
$$\alpha=\frac{2\mu_1\mu_2}{\mu_0+\mu_1+\mu_2-\mu_3-\mu_4},$$
and let
$$Q=(\mu_0X_1X_2+\mu_1X_0X_2+\mu_2X_0X_1)+\alpha X_0^2,$$
a singular quadric surface. Then as above, the intersection of $Q$
with $H$ consists of two lines, each counted twice, and a quartic
elliptic curve $F$. This quartic elliptic curve is the base locus of a
pencil of quadrics,
$$\langle Q,\alpha\mu_0X_3X_4+(\mu_1X_2+\mu_2X_1+\alpha
X_0)(\mu_3X_4+\mu_4X_3) \rangle.$$
Since we are looking, in the Kummer
case, for an element of this pencil that decomposes into two planes,
we begin by looking at the singular elements of the pencil. We find,
inter alia, the following result.

\begin{prop}
Let $\mu_0\mu_1\mu_2\mu_3\mu_4\not=0$, and assume
$$\alpha=\frac{2\mu_1\mu_2}{\mu_0+\mu_1+\mu_2-\mu_3-\mu_4},\qquad
   \beta=\frac{2\mu_3\mu_4}{\mu_0+\mu_3+\mu_4-\mu_1-\mu_2}$$
are finite. Then let
\begin{multline*}
R=(\mu_1X_2+\mu_2X_1)(\mu_3X_4+\mu_4X_3)+\alpha(\mu_0X_3X_4+\mu_3X_0X_4+\mu_4X_0X_3)\\
+\beta(\mu_0X_1X_2+\mu_1X_0X_2+\mu_2X_0X_1)+\alpha\beta X_0^2.
\end{multline*}
Then $R$ is always singular, i.e., of rank$\leq3$, with singular point
$$[\mu_1+\mu_2-\mu_3-\mu_4:-\mu_1:-\mu_2:\mu_3:\mu_4].$$
Further, $R$ has rank$\leq2$, i.e., decomposes, exactly if
$$\sum_{i=0}^4\mu_i^3-\sum_{i\not=j}\mu_i^2\mu_j
+2\sum_{i\not=j\not=k}\mu_i\mu_j\mu_k=0.$$
\end{prop}
\begin{proof}
This can all be checked by computing the necessary determinants.
\end{proof}

Observe that $R$ is symmetric with respect to the notational involution
\begin{align*}
*_1&\leftrightarrow *_3,\\
*_2&\leftrightarrow *_4,\\
\alpha&\leftrightarrow\beta,
\end{align*}
and in the case that $H$ is Kummer, we know that $R$ breaks into the
planes containing $C_{(03214)}$ and $C_{(01432)}$. Also, the
proposition provides for us a point on the intersection of these
planes. If let $S_5$ act on this proposition, so to speak, by
relabelling the variables, we obtain points on each of the
intersections of the six planes we are interested in. For example, the
plane $P_{(03214)}$ contains the following five points:
\begin{align*}
&[\mu_1+\mu_2-\mu_3-\mu_4:-\mu_1:-\mu_2:\mu_3:\mu_4]\\
&[-\mu_0:\mu_0+\mu_3-\mu_2-\mu_4:\mu_2:-\mu_3:\mu_4]\\
&[\mu_0:-\mu_1:\mu_2:\mu_1+\mu_4-\mu_0-\mu_2:-\mu_4]\\
&[\mu_0:\mu_1:-\mu_2:-\mu_3:\mu_2+\mu_3-\mu_0-\mu_1]\\
&[-\mu_0:\mu_1:\mu_0+\mu_4-\mu_1-\mu_3:\mu_3:-\mu_4]
\end{align*}
Taking minors of this matrix, we obtain an equation for the plane
$P_{(03214)}$ with coefficients cubic in the $\mu_i$s, and likewise
for the other five planes. More interestingly, we can find the
intersections of the conics
$$C_{(01324)},C_{(03421)},C_{(01432)}$$
with the line $\ell_{01}$, that is to say, the locations of the points
$p_{e+f-2a},p_{f+b-2a},p_{c+d-2a}$ on the trope $\Theta_{c+d-a}$.
This gives the following theorem.

\begin{thm}
If $\mu_0\mu_1\mu_2\mu_3\mu_4\not=0$ and
$$\sum_{i=0}^4\mu_i^3-\sum_{i\not=j}\mu_i^2\mu_j
+2\sum_{i\not=j\not=k}\mu_i\mu_j\mu_k=0,$$
and if the Hessian quartic surface $H$ given by
$$\V(\sum_{i=0}^4\mu_i\prod_{j\not=i}X_j)$$
has only ten nodes, then $H$ is Kummer.
Specifically, let $B$ be the branched cover of $\P^1$ over
\begin{align*}
a&=\frac{\mu_1+\mu_4-\mu_0-\mu_2-\mu_3}{2\mu_3},\\
b&=\frac{2\mu_2}{\mu_0+\mu_4-\mu_1-\mu_2-\mu_3},\\
c&=0,\\
d&=-1,\\
e&=\frac{\mu_0+\mu_3-\mu_1-\mu_2-\mu_4}{\mu_1+\mu_2-\mu_0-\mu_3-\mu_4},\\
f&=\infty.
\end{align*}
Then these six points will be distinct, so $B$ will be a smooth genus
2 curve. Let $K$ be its Kummer surface, and let $W$ be the Weber hexad
$$\{0,b+c-2a,c+d-2a,d+e-2a,e+f-2a,f+b-2a\}\subset K.$$
Then $H\cong\K$.

Conversely, if $\{a,b,c=0,d=-1,e,f=\infty\}$ are six distinct points
on $\P^1$, and $B$ is the genus 2 curve branched over those six
points, and $K$ and $W$ are as usual, the surface \K\ can be embedded
as a Hessian, with equation
$$H=\V(\mu_0X_1X_2X_3X_4+\mu_1X_0X_2X_3X_4+\mu_2X_0X_1X_3X_4+
\mu_3X_0X_1X_2X_4+\mu_4X_0X_1X_2X_3),$$
where the coefficients $\mu_i$ are given by
\begin{align*}
\mu_0&=a(b+1),&
\mu_1&=e(a+1),&
\mu_2&=b(a-e),&
\mu_3&=e-b,&
\mu_4&=(a-b)(e+1).
\end{align*}
\end{thm}
\begin{proof}
Given the previous theorem and the proposition, this reduces to a
computation.
\end{proof}

\section{Suggestions for further research}
As stated in the introduction, this exploration is by no means done.
To begin with, there is the problem of finding 6 points in $\P^2$ to
blow up to obtain the cubic surfaces associated to these Hessians.
Igor Dolgachev has presented a candidate sextuple, but this author
cannot see a good technique to answer his question.

Another intriguing line of research is broached by observing that
among all cubic surfaces, there is a codimension one subfamily of
singular cubic surfaces. On the moduli space, this divisor meets the
Kummer divisor studied in this paper, and the two divisors are
everywhere tangent along their intersection. However, the condition of
smoothness of the cubic surface has only barely made its presence
known in the results of this paper. A related question is brought up
by asking what happens if we allow our genus 2 curve to degenerate.

\section*{Acknowledgments}

We thank Igor Dolgachev for suggesting this problem. Also, we
gratefully acknowledge Bert van Geemen, who has been approaching the
same questions using theta function techniques. He independently found
the cubic relation described in this paper, and made several helpful
suggestions and comments while I was pursuing the geometric approach
to the problem.

\bibliographystyle{alpha} 
\bibliography{joelr}

\end{document}